\documentclass[11pt]{amsart} 
\title[Spectral synthesis for o. s. projective tensor product]{Spectral synthesis for operator space projective tensor product of $C^*$-algebras}

\usepackage[all]{xy}
\usepackage{amsmath,amsfonts,amssymb,amsthm,a4,hyperref}
\usepackage{enumerate}
\newtheorem{thm}{\sc Theorem}[section]
\newtheorem{cor}[thm]{\sc Corollary}
\newtheorem{prop}[thm]{\sc Proposition}
\newtheorem{lem}[thm]{\sc Lemma}
\theoremstyle{remark}
\newtheorem{defn}[thm]{\sc Definition}
\newtheorem{rem}[thm]{\sc Remark}
\newenvironment{pf}{\noindent {\sc Proof:}}{\hfill $\Box$}

\newcommand{\oop}{\widehat\otimes}
\newcommand{\seq}{\subseteq}
\newcommand{\oh}{\otimes^h}
\newcommand{\omin}{\otimes^{\min}}

\newcommand{\ra}{\rightarrow}
\newcommand{\ot}{\otimes}
\newcommand{\C}{\mathbb{C}}

\newcommand{\B}{\mathcal{B}}
\newcommand{\K}{\mathcal{K}}

\newcommand{\N}{\mathbb{N}}
\newcommand{\M}{\mathbb{M}}
\newcommand{\R}{\mathbb{R}}
\begin{document}

\author[R. Jain]{Ranjana Jain} \address{Department of Mathematics\\
  Lady Shri Ram College for Women\\ New Delhi-110024, India.}
\email{ranjanaj\_81@rediffmail.com}

\author[A. Kumar]{Ajay Kumar}
\address{Department of Mathematics\\ University of Delhi\\ Delhi-110007\\
  India.}  \email{akumar@maths.du.ac.in}
% ------------------------------------------------------------------------



\begin{abstract}
 We study the spectral synthesis for the Banach $*$-algebra $A\oop B$, the operator space projective tensor product of $C^*$-algebras $A$ and $B$.  It is shown  that if $A$ or $B$ has finitely many closed ideals, then $A\oop B$ obeys spectral synthesis.   The Banach algebra $A \oop A$ with the reverse involution is also  studied.
\end{abstract}

\maketitle 

\section{Introduction and notations}
For operator spaces $V$ and $W$, and $u\in  V\otimes W$, the {\it operator space projective tensor norm} \index{tensor norm/tensor product!operator space projective} is defined as
$$\|u\|_{\wedge}= \inf \{ \|\alpha\| \|v\| \|w\| \|\beta\| : u = \alpha (v\otimes w) \beta  \}, $$
where $\alpha \in \M_{1,pq}, \beta \in \M_{pq,1}, v\in M_p(V)$ and $w\in M_q(W)$, $p,q\in \N$ being arbitrary, and $v\otimes w = (v_{ij} \otimes w_{kl})_{(i,k),(j,l)} \in M_{pq}(V\ot W) $. The {\it operator space projective tensor product} $V\oop W$ is the completion of $V\otimes W$ under $\|\cdot\|_{\wedge}$-norm. The algebraic tensor product $V\otimes W $ is complete with respect to $\|\cdot\|_\wedge$-norm if and only if either $V$ or $W$ is finite dimensional. To see this,  if $V$ is finite dimensional, then as a Banach space it is isomorphic to $\C^n$, thus $V\oop W$ is Banach space isomorphic to the direct sum of n-copies of $W$, which is complete.  Also, if $V$ and $W$ are both infinite dimensional, then one can choose two sequences $\{e_n\}$ and $\{f_n\}$ of linearly independent vectors in $V$ and $W$ such that $\|e_n\|=\|f_n\|=1$ for all $n\in \N$. The sequence $(u_n)$ in $V\ot W$ defined as $u_n = \sum_{i=1}^n 2^{-i}e_i\otimes f_i$ is a Cauchy sequence with respect to $\|\cdot\|_\wedge$-norm, but is not convergent in $V\ot W$. It is known that for $C^*$-algebras $A$ and $B$, $A\oop B$ is a Banach $*$-algebra under natural involution (\cite{kumar}). 

The notion of spectral synthesis has been studied extensively for commutative and unital Banach algebras, for $L^1$-group algebras and for Banach $*$-algebras \cite{somer,somer2,somer3,kan}. Spectral synthesis for Banach space projective tensor product of commutative Banach algebras  and for the Haagerup tensor product of $C^*$-algebras has also been explored(\cite{kan,graham,ass,somer3}). Roughly speaking spectral synthesis holds for a Banach $*$-algebra $X$ if every closed ideal of $X$ is the intersection of primitive ideals containing it. Spectral synthesis for Banach space projective tensor product of commutative Banach algebras has already been explored(\cite{kan}). For commutative $C^*$-algebras $A$ and $B$, the natural contractive homomorphism of $A\oop B$ into $A\oh B$ is an isomorphism whose inverse has norm equal to Grothendieck constant. Thus, for countable locally compact Hausdorff spaces $X$ and $Y$, $C_0(X) \oop C_0(Y)$ has spectral synthesis. However, for cantor set  or any infinite compact group $D$, $C(D) \oop C(D)$ does not have spectral synthesis (\cite[11.2.1]{graham},\cite{kan}). 

In Section 2, we define the concept of spectral ideals in $A\oop B$, and prove that the Banach $*$- algebra $A\oop B$ has spectral synthesis if and only if each closed ideal of $A\oop B$ is spectral. This result is then used to produce plenty of spectral ideals in $A\oop B$.  We also discuss few cases where $A\oop B$ obeys spectral synthesis. In particular, we prove that if $A$ or $B$ has finitely many closed ideals, then $A\oop B$ has spectral synthesis. Thus,  the Banach $*$-algebras like  $C_0(X) \oop \B(H)$, $\B(H) \oop \K(H)$ and $\B(H) \oop \B(H)$ all obey spectral synthesis, $X$ being a locally compact topological space and $H$ being an infinite dimensional separable Hilbert space. In section 3, the algebra $A\oop B$ with the reverse involution is discussed. It is shown that with this involution the algebra is symmetric and $*$-semisimple only in the trivial cases.  

%Through out the discussion, $A$ and $B$ will denote $C^*$-algebras, until unless specified. 
For a Banach algebra $X$, we denote the set of closed (two-sided) ideals of $X$ by $Id(X)$,  the set of proper closed ideals of $X$ by $Id^\prime(X)$ and the set of all prime ideals by $Prime(X)$. If $X$ is a Banach $*$-algebra, then $Prim(X)$ stands for the set of primitive ideals of $X$, that is, the set of all kernels of irreducible $*$-representations of $X$ on Hilbert space. There is a topology $\tau_w$ on $Id(X)$ which is generated by the sub-basic open sets of the form 
$$ Z_J :=\{ I \in Id(X): I\nsupseteq J\},\, J \in Id(X).$$
We throughout use the notation $q_J$ for the quotient map $q_J: A \ra A/J$. Recall that, for closed ideals $M$ and $N$ of $C^*$-algebras $A$ and $B$, the map $q_M\ot q_N : A\ot B\ra A/M\ot B/N$ extends to quotient maps $q_M\oop q_N: A \oop B \ra A/M \oop B/N$ and $q_M\omin q_N: A \omin B \ra A/M \omin B/N$.

Let $A$ and $B$ be $C^*$-algebras. Define a map $\Phi: Id(A) \times Id(B) \ra Id(A\oop B)$ as
$$\Phi (M,N) = A\oop N + M \oop B. $$
The map $\Phi$ is well defined by \cite[Proposition 3.2]{ranj2}. It satisfies many nice topological properties listed as below:
\begin{prop}\label{spec5}
Let $A$ and $B$ be $C^*$-algebras and $\Phi: Id(A) \times Id(B) \ra Id(A\oop B)$ be defined as above. Then
\begin{enumerate}[(i)]
 \item\label{2} $\Phi$ maps $Prime(A)\times Prime(B)$ onto $Prime (A\oop B)$.
 \item $\Phi$ maps $Prim(A)\times Prim(B)$ into $Prim(A\oop B)$. If $A$ and $B$ are separable, then $\Phi$ maps $Prim(A)\times Prim(B)$ onto $Prim (A\oop B)$.
\item $\Phi$ maps $Id^\prime(A)\times Id^\prime(B) $ into $Id^\prime(A\oop B)$ injectively.
\item The mapping $\Phi$ is $\tau_w$-continuous.
\item The restriction of $\Phi$ to $Id^\prime (A) \times Id^\prime(B)$ is a homeomorphism onto its image in $Id^\prime(A\oop B)$.
\item  The restriction of $\Phi$ to $Prime (A) \times Prime(B)$ is a homeomorphism onto $Prime(A\oop B)$.
\end{enumerate}
\end{prop}

\begin{pf}
(i) and (ii) follow from Theorems 3.1 and 3.2 of \cite{ranj3}, respectively.

For (iii), note that, for proper closed ideals $M$ and $N$ of $A$ and $B$, the isomorphism of $A/M \oop B/N$ onto $(A\oop B)/(A\oop N + M\oop B)$ (\cite[Lemma 2.2]{ranj3}) assures that  $A\oop N + M\oop B$ is also proper in $A\oop B$. Further,  for $M_1,M_2 \in Id^\prime (A), N_1,N_2 \in Id^\prime (B)$, $A\oop N_1 + M_1\oop B \subseteq A\oop N_2 + M_2\oop B$ if and only if $M_1\seq M_2 , N_1 \seq N_2 $.  To see this, consider an $m \in M_1$, so that for an arbitrary $b\in B$, $m\ot b \in \ker(q_{M_1}\oop q_{N_1}) \seq \ker(q_{M_2}\oop q_{N_2})$, giving $q_{M_2}(m) = 0$, that is $m\in M_2$, and similarly $N_1 \seq N_2$. Thus, $\Phi$ is injective.

(iv)-(vi) can be proved exactly on the same lines of their counterparts in Haagerup tensor product as discussed in Lemma 1.4 and Theorem 1.5 of \cite{arch}.   
\end{pf}

\vspace*{1mm}
Throughout this paper $A$ and $B$ represent $C^*$-algebras, until otherwise specified.
%========================================================
\section{Spectral synthesis}

We first give the standard definition of spectral synthesis for a Banach $*$-algebra that appear in the literature.  Let $X$ be a Banach $*$-algebra. For each $E\subseteq Prim(X)$, there associates a closed ideal  {\it kernel} of $E$ defined as
$$ k(E) := \cap_{P \in E} P.$$
Also, for each  $M\subseteq X$, {\it hull} of $M$ is defined as
$$ h_X(M) := \{ P \in Prim(X) : P \supseteq M \}.$$
We shall denote the hull of $M$ by $h(M)$, when there is no confusion with $X$. Equip $Prim(X)$ with the {\it hull-kernel topology} (or, hk-topology), where for every $E\subseteq Prim(X)$, its closure is $\overline{E} = h(k(E))$. Similarly, one can talk about the hk-topology on $Prime(X)$.  Note that, if $E\subseteq Prime(X)$, then the relative $\tau_w$-topology on $E$ coincides with the hull-kernel topology.

\begin{defn}\label{defII}
 A closed subset $E$ of $Prim(X)$ is called {\it spectral} if $k(E)$ is the only closed ideal in $X$ with hull equal to $E$.  A Banach $*$-algebra $X$ is said to have {\it spectral synthesis} if every closed subset of $Prim(X)$ is spectral. 
\end{defn}

A closed ideal of Banach $*$-algebra $X$ is said to be semisimple if it is an intersection of all the primitive ideals of $X$ containing it.
\begin{prop}\label{prim}
Let $X$ be a Banach $*$-algebra having Wiener property. Then $X$ has spectral synthesis if and only if for every $J \in Id(X)$, $ J =k(h(J))$, or, in other words, every closed ideal of $X$ is semisimple.   
\end{prop}

\begin{pf}
Let us consider a proper closed ideal $J$ of $X$.  Since $X$ has Wiener property, there exists an irreducible $*$-representation, say $\pi$, of $X$ which annihilates $J$, that is, $J\subseteq \ker \pi$, so that $E = h(J)$ is non empty. We claim that $E$ is closed in the hk-topology.
Let $ Q\in \overline{E}=h(k(E))$, then $k(E) \seq Q$. Since $J \seq P$ for all $P\in E$ we have $J\seq k(E) \seq Q$, so that $Q\in E$. which gives that
%$\overline{E}\seq E$, giving that
$E$ is closed. Since $X$ obeys spectral synthesis, and $E = h(J)$, we have $J= k(E)$, that is, $J$ is the intersection of primitive ideals containing it. Also, note that since $X$ has Wiener property, the empty set $\phi$ is spectral, so that $X= k(h(X))$.
% $\phi spectral, meaning is $k(\phi) is the only closed ideal of $A$ with hull equal to $\phi$, but $h(A) = \phi$, so that $A= k(\phi)$. 

Converse follows easily from the given condition.
%For the converse, consider a closed subset $E$ of $PrimX$, that is, $E = h(k(E))$. Let $J$ be a closed ideal of $A\oop B$ with $h(J) =E $. By the given condition $J= k(h(J))= k(E)$. 
\end{pf}

\begin{cor}\label{prim1}
 Let $X$ be a Banach $*$-algebra having Wiener property. Then $X$ has spectral synthesis if and only if there is a one-one correspondence between the closed ideals of $X$ and the $\tau_w$-open subsets of $Prim(X)$ (or, $Prime(X)$). 
\end{cor}

\begin{pf}
Let $X$ have spectral synthesis. For $J\in Id(X)$, recall $Z_J:=\{ P\in Prim(X):  P\nsupseteq J\} = Prim(X) \setminus h(J)$ is an open subset of $Prim(X)$ under the relative $\tau_w$-topology, so that we have a well defined correspondence $J \mapsto Z_J$ between the closed ideals of $X$ and $\tau_w$-open subsets of $Prim(X)$.
 For $K, L \in Id(X)$, it is clear from Proposition \ref{prim} that $K = k(h(K))$, and $L = k(h(L))$. Thus, it can be easily seen that 
$$K \seq L \quad \text{if and only if} \quad Z_K \seq Z_L,$$
%The necessary part is clear from the definition. Conversely, by Remark \ref{spec4}(\ref{1}), we have $K = \cap \{ P\in X: P\supseteq K\} = k(h(K))$, and $L = k(h(L))$. Since $Z_K \seq Z_L$, by definition, $h(L) \seq h(K)$, which futher implies, $K \seq L$.
%for prime ideal $P$ of $A\oop B$, $K\nsubseteq P \Rightarrow L \nsubseteq P$, that is $L\subseteq P \Rightarrow K \subseteq P$. So, $\{ P: L\seq P\} \seq \{P: K\seq P\}$, giving that $K\seq L$. 
which shows that the correspondence in one-one. Now consider a $\tau_w$-open subset $G$ of $Prim(X)$, and set $J := k(Prim(X)\setminus G)$. Since $Prim(X)\setminus G$ is closed under the hull-kernel topology, 
 $$Z_J = Prim(X)\setminus h(k(Prim(X)\setminus G)) = Prim(X) \setminus (Prim(X)\setminus G) = G,$$ 
which proves that this correspondence is surjective.

Conversely, for every closed ideal $I$ of $X$, since $h(I) = h(k(h(I))$, we have $Z_I = Z_{k(h(I))}$. Using the given condition, this gives $I = k(h(I))$. Result now follows from Proposition \ref{prim}.
\end{pf}

\begin{rem}\label{spec4}
For $C^*$-algebras $A$ and $B$, since $A\oop B$ has Wiener property (\cite[Theorem 4.1]{ranj2}), $A\oop B$ has spectral synthesis if and only if every closed ideal $J$ of $A\oop B$ is semisimple. In particular,  if $A\oop B$ has spectral synthesis then every closed ideal $J$ of $A\oop B$ is the intersection of prime ideals containing $J$. 
% \item It is known that for a $C^*$-algebra $A$, the class $\{h(M):M\seq A\}$ forms the closed sets for the hull-kernel topology on $Prim(A)$. Also there is one-one correspondence between the open sets of $Prim(A)$ under the hull-kernel topology and the closed ideals of $A$  (\cite[Theorem 4.1.3]{ped}).  From the above result, one can notice that the both the statements are true for  for any Banach $*$-algebra $X$, the collection $\{ h(M) :M\seq X\}$ forms the closed sets for the hull-kernel topology of $Prim(X)$.
\end{rem}

%\vspace*{2mm}
The next two results connect the spectral synthesis of a Banach $*$-algebra with that of its ideal and the corresponding quotient algebra. The first result follows on the similar lines as that in \cite[Proposition 1.16]{somer3}. However, we present here a proof for the sake of completion. 

\begin{prop}\label{spec15}
 Let $X$ be a Banach $*$-algebra with Wiener property, and $J$ be a closed $*$-ideal of $X$ having bounded approximate identity and Wiener property. If $J$ and $X/J$ both have spectral synthesis (as Banach $*$-algebras), then $X$ has spectral synthesis. 
\end{prop}

\begin{pf}
By Corollary \ref{prim1}, it is sufficient to show that for $I, K \in Id(X)$, $I=K$, whenever $h_X(I)= h_X(K)$. Note that, since $X$ has Wiener property, $X/J$ also has Wiener property, so by Proposition \ref{prim}, every  closed ideal of $J$ and $X/J$ is semisimple.  For $P \in h_{X/J}(q_J(I))$, $P= \ker \pi$ with $\pi(q_J(I)) = \{0\}$, $\pi:X/J \ra \B(H)$ being an irreducible $*$-representation. Then  $\pi_0:= \pi \circ q_J$ is an irreducible $*$-representation of $X$ on $H$ with $\pi_0(I)=0$. Since $h_X(I)= h_X(K)$, $\pi_0 \in h_X(K)$, which further gives $P\in h_{X/J}(q_J(K))$. Thus, $h_{X/J}(q_J(I))= h_{X/J}(q_J(K))$. Since $J$ has an approximate identity, by \cite[Proposition 2.4]{dixo}, $I+J$ and $K+J$ are closed in $X$, so that $q_J(I)$ and $q_J(K)$ are closed ideals of $X/J$. Since $X/J$ obeys spectral synthesis, by Proposition \ref{prim}, $q_J(I)= q_J(K)$. Further,  for any closed ideal $L$ of $J$, it is routine to check that there is a one-one correspondence between the sets $\{P \in h_X(L): J \nsubseteq P\} $ and $h_J(L)$ via $P \mapsto P\cap J$. Xlso, $h_X(I\cap J) = h_X(I) \cup h_X(J)= h_X(K) \cup h_X(J)= h_X(K\cap J)$. Thus, it can be easily seen that  $h_J(I\cap J) = h_J(K\cap J)$.  Since $J$ has spectral synthesis, this gives, $I\cap J = K\cap J$.

Now, consider $x\in I$, then $q_J(x) = q_J(y)$ for some $y \in K$, so that $a:=x-y \in J$. Let $J_a$ be the smallest closed ideal of $J$ containing $a$. Since $J$ obeys spectral synthesis, $J_a=\cap \{P\in Prim(J): J_a \seq P\}$. Clearly $\overline{JaJ}\seq J_a$. Now consider $P \in Prim(J)$ such that $JaJ \seq P$. Since $P$ is prime being primitive, this gives $a\in P$ which shows that $J_a \seq P$. Thus $\overline{JaJ} = \cap \{P\in Prim(J): JaJ \seq P\} = \cap \{P\in Prim(J): J_a \seq P\} = J_a $. So
\begin{eqnarray*}
 x-y \in \overline{J(x-y)J} & \seq & \overline{JxJ -JyJ} \\
           & \seq & \overline{JIJ - JKJ}\\
        & \seq & \overline{I\cap J - K\cap J} \\
       & =&  K\cap J. 
\end{eqnarray*}
So, $x = y-(y-x) \in K+ (K\cap J) =K$, which gives $I\seq K$. Similarly, $K\seq I$, which proves the claim. 
\end{pf}

\vspace*{2mm}
In fact, the converse of the above statement is also true as presented below.
\begin{prop}
Let $X$ be a Banach $*$-algebra with a closed $*$-ideal $J$ such that $X$ and $J$ both possess Wiener property. If $X$ obeys spectral synthesis, then so does $J$ and $X/J$.  
\end{prop}

\begin{pf}
 By Proposition \ref{prim}, it is enough to check that for a closed ideal $L$ of $J$, $L = k(h_J(L))$. Since every closed ideal of $X$ is semisimple, and every primitive ideal is prime, from \cite[Proposition 1.14]{somer3}, $L$ is also a closed ideal of $X$, so that by Proposition \ref{prim},  $L = k(h_X(L)) $.    It can be easily verified that there is a one-one correspondence between the sets $\{P \in h_X(L): J \nsubseteq P\} $ and $h_J(L)$ via $P \mapsto P\cap J$. So, we have
\begin{eqnarray*}
 L = L \cap J & = & \bigcap_{P \in h_X(L)} (P\cap J) \\ 
        & = & \bigg(\bigcap_{\substack{P \in h_X(L) \\ J \nsubseteq P}} (P\cap J)\bigg) \cap \bigg(\bigcap_{\substack{P \in h_X(L) \\ J \subseteq P}} (P\cap J)\bigg)\\
  & = & \big(\bigcap_{P^\prime \in h_J(L)} P^\prime\big) \cap J \\
  & = & k(h_J(L))
\end{eqnarray*}
Thus, $J$ obeys spectral synthesis.

Next, consider a closed ideal $K$ of $X/J$. Since $X/J$ has Wiener property, it is enough to check that $K \supseteq k(h_{X/J}(K))$. Consider an element $x \in  k(h_{X/J}(K))$, where $x= y+J \in X/J$. Note that $K= I/J$ for some closed ideal $I$ of $X$ containing $J$. Using the one-one correspondence between $Prim(X/J)$ and   $\{ P \in Prim(X): J\seq P \}$, one can check that $y\in  k(h_X(I))$. Since $X$ has spectral synthesis, $I = k(h_X(I))$, so that $y \in I$, which shows that $x \in K$. Hence the result.
% For any $P \in h_{A/J}(K)$
\end{pf}

\vspace*{2mm}
We are now prepared to discuss spectral synthesis for operator space projective tensor product $A\oop B$ of $C^*$-algebras $A$ and $B$.
Allen, Sinclair and Smith, in \cite{ass}, defined the concept of spectral synthesis for the Haagerup tensor product of $C^*$-algebras in a somewhat different flavor. In the same spirit, using the terminologies of \cite{ass}, we give another definition for the spectral synthesis of $A\oop B$. It is known that for any $C^*$-algebras $A$ and $B$, the canonical $*$-homomorphism $i: A\widehat{\otimes}B \ra A \omin B$ is injective (\cite[Corollary 1]{ranj1}), so that we can regard $A\oop B$ as a $*$-subalgebra of $A\omin B$. Consider a closed ideal $J$ of $A\oop B$ and let $J_{\min}$ be the closure of $i(J)$ in $A\omin B$, in other words, $J_{\min}$ is the min-closure of $J$ in $A \omin B$. Now we associate two closed ideals, namely the  upper and the lower ideals, with $J$ as:
%As done in \cite{ass}, given a closed ideal $K$ of $A \omin B$, there associate with $K$ two closed ideals $K_l$ and $K^u$ in $A\oop B$ defined by 
\begin{eqnarray*}
 J_l & = & \, \text{closure of span of all elementary tensors of} \, J \,\, \text{in} \, A\oop B, \\
J^u & = & J_{\min} \cap (A\oop B).
\end{eqnarray*}
%$$J_l = (J_{\min})_l, J^u = (J_{\min})^u$$
Clearly $J_l \seq J \seq J^u$ for any closed ideal $J$ of $A\oop B$.
\begin{defn}\label{defI}
  A closed ideal $J$ of $A\widehat{\otimes}B$ is said to be  {\it spectral} if $J_l=J=J^u$. 
%The algebra $A\oop B$ is said to have {\it spectral synthesis}\index{spectral synthesis!for $A\oop B$} if every closed ideal of $A\oop B$ obeys spectral synthesis.
\end{defn}

The main aim of this section is to show that $A\oop B$ has spectral synthesis if and only if its every closed ideal is spectral. We first characterize the upper ideals in terms of primitive ideals. 
\begin{lem}\label{spec0}
For closed ideals $M$ and $N$ of $A$ and $B$,
$$\ker(q_M \oop q_N) = \ker(q_M \omin q_N)\cap A\oop B.$$
\end{lem}

\begin{pf}
For $z\in A\oop B$, let $\{z_n\}$ be a sequence in $A\ot B$ such that $\lim_n \|z_n -z\|_\wedge =0$, then
\[  \|(q_M \oop q_N)(z_n) - (q_M \oop q_N)(z)\|_{\min} \leq \|q_M \oop q_N\| \|z_n -z\|_\wedge, \]
which shows that the sequence $ \{(q_M \oop q_N)(z_n)\} $ is convergent to  $(q_M \oop q_N)(z)$ in $A\omin B$.
 Also, $\|z_n -z\|_{\min} \leq \|z_n -z\|_\wedge$, so that $\lim_n \|z_n -z\|_{\min} =0$, which further gives 
$$ (q_M \omin q_N)(z_n) \stackrel{\min}{\longrightarrow} (q_M \omin q_N)(z).$$
Since both the mappings $q_M \oop q_N$ and $q_M \omin q_N$ agree on $A\ot B$, by continuity, we have $(q_M \omin q_N)(z) = (q_M \oop q_N)(z)$, and this is true for all $z\in A\oop B$, proving the given relation.
\end{pf}

\begin{prop}\label{ss1}
 For a closed ideal $J$ of $A\oop B$, $J = J^u$ if and only if $J$ is semisimple. 
\end{prop}

\begin{pf}
Let us first assume that $J = J^u$. Since, in a $C^*$-algebra every closed ideal is semisimple,  $J_{\min} = \cap \{\ker \tilde{\pi}_\alpha: J_{\min} \seq \ker \tilde{\pi}_\alpha \}$, where each  $\tilde{\pi}_\alpha$ is an irreducible $*$-representation of $A\omin B$ on some Hilbert space. Set $\pi_\alpha :=  \tilde{\pi}_\alpha \circ i$, then each $\pi_\alpha$ is an irreducible $*$-representation of $A\oop B$ annihilating $J$. Using some routine  calculations, and the fact that $J= J^u$ one can prove that $J = \cap \ker \pi_\alpha$. Note that, although the collection $\{P 
\in Prim(A\oop B): J\seq P\}$ is larger than $\{ \ker \pi_\alpha : \pi_\alpha =  \tilde{\pi}_\alpha \circ i\}$, it is easy to check that  $J$ is actually the intersection of all the primitive ideals of $A\oop B$ containing $J$.
%Let $S_1= \{P \in Prim(A\oop B): J\seq P\}$ and $S_2 = \{ \ker \pi_\alpha : \pi_\alpha =  \tilde{\pi}_\alpha \circ i\}$. Then, clearly $S_2 \seq S_1$ so that $\cap_{P \in S_1} P \seq \cap_{P \in S_2 } P = J$. Also for any $P \in S_1, J \seq P$ so that $J \seq \cap_{P_\in S_1} P$. Thus $J = \cap_{P \in S_1} P$. 

Conversely, let $J= \cap_{J\seq P_\alpha} P_\alpha$, $P_\alpha $ being primitive ideals of $A\oop B$.  Let, if possible, there exist an element  $x\in J^u $ such that $x \notin J$. Then $x\notin P_\alpha$ for some $\alpha$. Since $P_\alpha$ is primitive, by \cite[Theorem 3.2]{ranj3}, there exist closed (prime) ideals $M$ and $N$ in $A$ and $B$, respectively, such that $P_\alpha = A\oop N + M\oop B$. Now, consider the bounded homomorphisms $q_M \oop q_N :A\oop B \ra A/M \oop B/N$, and $q_M \omin q_N :A\omin B \ra A/M \omin B/N$ with $\ker(q_M \oop q_N) = P_\alpha$ (\cite[Proposition 3.5]{ranj2}. By Lemma \ref{spec0}, $x\notin \ker(q_M \omin q_N)$, which by Hahn Banach Theorem gives a $\phi \in (A\omin B)^*$ such that $\phi(x) \neq 0$, and $\phi(\ker(q_M \omin q_N)) =\{0\}$. The relation  $J \seq P_\alpha \seq \ker(q_M \omin q_N)$ gives $J_{\min} \seq \ker(q_M \omin q_N)$, which further shows that $\phi(J_{\min}) =0$. Thus $x\notin J_{\min}$, which gives a contradiction to the fact that $x\in J^u$. Hence the result.
\end{pf}

\vspace*{2mm}
Using Propositions \ref{prim} and \ref{ss1}, we have a following characterization for spectral synthesis in terms of upper ideals. 
\begin{thm}\label{spec2}
 The Banach $*$-algebra $A\oop B$ has spectral synthesis if and only if $J= J^u$, for every closed ideal $J$ of $A\oop B$.
\end{thm}
We now prove that the Banach $*$-algebra $A\oop B$ has spectral synthesis if and only if every closed ideal of $A\oop B$ is spectral. We borrow some ideas from \cite{lazar} to prove the same. We first need an elementary result.
\begin{lem}\label{spec1}
 Let $J_i$ and  $K_i$ be closed ideals of $C^*$-algebras $A_i$, $i=1,2$. Then $J_1\oop J_2 \seq A_1\oop K_2 + K_1 \oop A_2$ if and only if either $J_1 \seq K_1$ or $J_2\seq K_2$.  
\end{lem}

%\begin{pf}
%For the necessary part, let, if possible, neither $J_1 \nsubseteq K_1$ nor $J_2\nsubseteq K_2$. Then, there exist $j_1\in J_1$, $j_2\in J_2$ such that $q_{K_1}(j_1) \neq 0, q_{K_2}(j_2) \neq 0$. This gives $(q_{K_1}\oop q_{K_2})(j_1\ot j_2) \neq 0$, thus $j_1 \ot j_2\notin \ker (q_{K_1}\oop q_{K_2})$, giving that $J_1\oop J_2 \nsubseteq A_1\oop K_2 + K_1 \oop A_2$. Converse is trivial.
%\end{pf}

\begin{thm}\label{spec6}
 For $C^*$-algebras $A$ and $B$, the Banach $*$-algebra $A\oop B$ has spectral synthesis if and only if every closed ideal of $A\oop B$ is spectral.
\end{thm}

\begin{pf}
We just need to prove that for every closed ideal $J$ of $A\oop B$,  $J = J_l$, if $A\oop B$ has spectral synthesis. 
Using Corollary \ref{prim1}, it is sufficient to show that $Z_J \seq Z_{J_l}$, where $Z_J:=\{ P\in Prime(A\oop B):  P\nsupseteq J\}$. Set $X:= Prime(A\oop B)$ and consider an element $P$ of $ Z_J$. Since $Z_J$ is an open subset of $X$ and $\Phi:Prime(A) \times Prime(B) \ra X$ is continuous, there exist open subsets $U_1$, $U_2$ of $Prime(A)$ and $Prime(B)$ such that $\Phi(U_1 \times U_2) \seq Z_J$ and $P\in \Phi(U_1\times U_2)$. Let $J_1 \in Id(A)$, $J_2\in Id(B) $  be the corresponding closed ideals such that $U_i= Z_{J_i}$, $i=1,2$.
We claim that $Z_{J_1\oop J_2} = \Phi(U_1 \times U_2)= \Phi(Z_{J_1} \times Z_{J_2})$. For any $Q \in Z_{J_1\oop J_2}$, by definition, $J_1\oop J_2 \nsubseteq Q$. Since $Q\in X$, and $\Phi$ is onto (Proposition \ref{spec5}), there exists $Q_1 \in Prime(A)$, $Q_2\in Prime(B)$ such that $A\oop Q_2 + Q_1 \oop B = \Phi (Q_1, Q_2) = Q $. By Lemma \ref{spec1}, $J_1\nsubseteq Q_1$ and $J_2 \nsubseteq Q_2$. This implies that $Q_i \in Z_{J_i} = U_i$, so that $Q=\Phi(Q_1 , Q_2) \in \Phi(U_1\times U_2)$. Thus,   $Z_{J_1\oop J_2} \seq \Phi(U_1 \times U_2)$. For the other containment, consider $\Phi(K_1 , K_2) \in \Phi(Z_{J_1} \times Z_{J_2})$. Since $K_i \in Z_{J_i}$, we have $J_i\nsubseteq K_i$, so by Lemma \ref{spec1}, $J_1\oop J_2 \nsubseteq \Phi(K_1,K_2)$. Note that $\Phi(K_1,K_2) \in X$, thus by the definition, $\Phi(K_1,K_2) \in Z_{J_1\oop J_2}$. 
So, $Z_{J_1\oop J_2} \seq Z_J$, which further gives, $ J_1\oop J_2 \seq J$. But, the definition of $J_l$ says that $J_1\oop J_2 \seq J_l$. This means that $Z_{J_1\oop J_2} \seq Z_{J_l}$. Since, $P \in \Phi(U_1\times U_2) = Z_{J_1\oop J_2}$, this gives $P\in Z_{J_l}$. Thus $Z_J \seq Z_{J_l}$, which proves that $J \seq J_l$, and hence the result. 
\end{pf}

\begin{rem}\label{spec16}
In other words, if $A\oop B$ obeys spectral synthesis, then every closed ideal $J$ of $A\oop B$ is the closure of the sum of all product ideals $J_1 \oop J_2 \seq J$, where $J_1 \in Id(A)$, $J_2 \in Id(B).$     
\end{rem}

The Banach $*$-algebra $A\oop B$ contains plenty of spectral ideals as demonstrated in the following and some later examples.
\begin{prop}\label{spec7}
For $I\in Id(A)$ and $J \in Id(B)$, the closed ideal $A \oop J+ I\oop B$ of $A\oop B$ is spectral. In particular, every closed maximal ideal, primitive ideal and prime ideal of $A\oop B$ is spectral.
\end{prop}

\begin{pf}
Set $K := A\oop J + I \oop B = \ker(q_I\oop q_J)$, then it is clear from the definition that $K = K_l$. Consider an element $ u\in K^u$. Let, if possible, $u\notin K$, then by Lemma \ref{spec0}, $u \notin \ker(q_I\omin q_J)$. Now, $K \seq \ker(q_I\omin q_J)$ implies $K_{\min} \seq \ker(q_I\omin q_J)$, giving $ u \notin K_{\min}$, a contradiction. Thus $K$ is spectral. Rest follows from the fact that every maximal, primitive and prime ideal can be expressed as an ideal of this form (\cite[Theorem 3.10]{ranj2}, \cite[Theorem 3.1, 3.2]{ranj3}).  
\end{pf}

Next, we prepare the ingredients to prove that for an infinite dimensional separable Hilbert space $H$, $\B(H) \oop \B(H)$ obeys spectral synthesis. We first need some elementary results regarding the lower and upper ideals of a closed ideal. 
\begin{prop}
For closed ideals $J$ and $K$ in $A \widehat{\otimes} B$, we have:
\begin{enumerate}[(a)]
\item $J_l \subseteq K_l$ and $J^u \subseteq K^u$, if $J \subseteq K$;
\item $(JK)_l= J_lK_l =J_l \cap K_l = (J\cap K)_l$, if $J_l$ or $K_l$ has a bounded  approximate identity;
\item $(J \cap K )^u \subseteq J^u \cap K^u$, with equality if $J = J^u, K = K^u$.
\end{enumerate}
\end{prop}

\begin{pf} (a) is trivial.  For (b), we first show that $J_l \cap K_l \seq J_lK_l$. Let $x\in J_l \cap K_l$ and assume that $J_l$ has bounded approximate identity. By Cohen's Factorization Theorem, there exist $y,z \in J_l$ such that $x= yz $ and $z$ belongs to the closed left ideal generated by $x$ in $J_l$. Clearly, $z \in J_l \cap K_l$, so that $x\in J_lK_l$. Thus, $J_l\cap K_l \seq J_lK_l$. Now, for an elementary tensor $x$ in $J\cap K$, clearly $x \in J_l \cap K_l$, giving $(J\cap K)_l \seq J_l\cap K_l$.
Also, for $a=\sum_{i=1}^n x_i\ot y_i \in J_l$ and $b=\sum_{j=1}^m z_j\ot w_j \in K_l$, clearly $ab \in (JK)_l$, being an elementary tensor of $JK$. Since $J_l$ and $K_l$ are both generated by elementary tensors, routine calculations show that $J_lK_l \seq (JK)_l$. Thus, we have
$$ (J\cap K)_l \seq J_l\cap K_l \seq J_lK_l \seq (JK)_l \seq (J\cap K)_l,$$
which gives the required equality. For (c), using the fact that $(J\cap K)_{\min} \seq J_{\min} \cap K_{\min} $, we get
\[
  (J \cap K)^u \seq J_{\min} \cap K_{\min} \cap A\oop B = J^u \cap K^u. 
\] 
\end{pf}

Following are  some direct consequences of the above proposition. 
\begin{cor}\label{spec8}
 If $I$ and $J$ are closed ideals of $A\oop B$ with at least one of them having bounded approximate identity. Then $I\cap J$ is spectral, whenever $I$ and $J$ are spectral.
\end{cor}

\begin{cor}\label{spec9}
 Every product ideal of $A\oop B$ is spectral. In particular, for closed ideals $I$ and $J$ of $A$ and $B$, $I\oop J = (I\omin J) \cap A\oop B$.
\end{cor}

\begin{pf}
For a product ideal $I \oop J$ of $A\oop B$, using \cite[Proposition 2.4]{ranj3}, we can write $$I\oop J  = (A\oop J) \cap (I\oop B).$$ Also, from \cite[Lemma 3.1]{ranj2}, $A\oop J$ and $I\oop B$ both possess bounded approximate identities. Thus, from Proposition \ref{spec7} and Corollary \ref{spec8}, $I\oop J$ is spectral.  
Clearly, 
$$I\oop J = (I\oop J)^u = (I \oop J)_{\min} \cap A\oop B = (I\omin J ) \cap A\oop B.$$ 
\end{pf}

\begin{cor}\label{spec13}
If either $A$ or $B$ is a simple $C^*$-algebra, then $A\oop B$ obeys spectral synthesis.
\end{cor}

\begin{pf}
 Let $A$ be simple. By \cite[Theorem 3.8]{ranj2}, every closed ideal of $A\oop B$ is a product ideal and thus is spectral by Corollary \ref{spec9}. Using Theorem \ref{spec6}, we get $A\oop B$ obeys spectral synthesis.
\end{pf}

In particular, for any $C^*$-algebra $A$, the Banach $*$-algebras $A\oop C_r^*(\mathbb{F}_2)$, $A\oop A_\infty$  and $A\oop \K(H)$ obey spectral synthesis, where  $C_r^*(\mathbb{F}_2)$ is the $C^*$-algebra associated to the left regular representations of the free group $\mathbb{F}_2$ on two generators, $A_\infty$ is the Glimm algebra (\cite{ped}) and $\K(H)$ is the $C^*$-algebra of compact operators on an infinite dimensional separable Hilbert space $H$.

\begin{thm}
 For an infinite dimensional separable Hilbert space $H$, the Banach $*$-algebra $\B(H) \oop \B(H)$ obeys spectral synthesis.
\end{thm}

\begin{pf}
 From \cite[Theorem 3.11]{ranj2}, we know that the only non trivial closed ideals of $\B(H) \oop \B(H)$ are $\K(H)\oop \K(H)$, $\B(H) \oop \K(H), \K(H) \oop \B(H)$ and $\B(H) \oop \\ \K(H) + \K(H) \oop \B(H)$. Using Proposition \ref{spec7} and Corollary \ref{spec9}, we can see that all the proper closed ideals of $\B(H) \oop \B(H)$ are spectral. The result now follows from Theorem \ref{spec6}.
\end{pf}

\begin{prop}\label{spec14}
Let $A$ and $B$ be $C^*$-algebras such that $A$ or $B$  has finitely many closed ideals. Then $A\oop B$ obeys spectral synthesis.
\end{prop}  

\begin{pf}
Without loss of generality, we may assume that $B$ has finitely many closed ideals say $n$, where $n\geq 2$. We prove the result by induction on $n$. For $n=2$, $B$ is simple and the result follows from Corollary \ref{spec13}. Let the result be true for all $C^*$-algebras with at most $(n-1)$ ideals. Let $B$ have $n >2$ closed ideals. Since there are finitely many closed ideals of $B$, there exists a minimal (non-trivial) closed ideal, say $K$, of $B$, which is clearly simple. Consider the closed $*$-ideal $J:=A\oop K$ of $X:=A\oop B$. Since $K$ is simple, using Corollary \ref{spec13}, it is clear that $J$ has spectral synthesis. Note that, by  \cite[Lemma 2.2(1)]{ranj3}, $X/J$ is isomorphic to $A\oop (B/ K)$ and the latter has spectral synthesis by induction hypothesis, since $B/K$ has atmost $(n-1)$ closed ideal. So, $X/J$ also has spectral synthesis. Moreover, $J$ and $X/J$ both have Wiener property (\cite[Theorem 4.1]{ranj3}), and $J$ has bounded approximate identity (\cite[Lemma 3.1]{ranj2}), the result now follows from Proposition \ref{spec15}.     
\end{pf}

\vspace*{2mm}
Thus, for any $C^*$-algebra $A$, $A\oop \B(H)$ obey spectral synthesis, where $H$ is a separable infinite dimensional Hilbert space. In particular, $C_0(X) \oop \B(H)$, $\B(H)\oop \B(H)$ and $\B(H)\oop \K(H)$ obey spectral synthesis, where $X$ is a locally compact Hausdorff space.  For more examples of  $C^*$-algebras with finitely many closed ideals, see \cite{lin}.
%The following consequence of the above result provides alternate way to obtain all the closed ideals of $\B(H)\oop \B(H)$, which were derived in \cite{ranj2} by another method.
\begin{cor}
 If $A$ and $B$ both have finite number of closed ideals, then every closed ideal of $A\oop B$ is a finite sum of product ideals. 
\end{cor}

\begin{pf}
 It follows from Proposition \ref{spec14} and Remark \ref{spec16}.
\end{pf}

%\begin{prop}\label{spec11}
%Let $A$ and $B$ be $C^*$-algebras such that $A \otimes B$ is $*$-regular. Then $A\oop B$ is $*$-regular. In particular, $A\oop B$ is $*$-regular, whenever $A$ is nuclear. 
%\end{prop}

%\begin{pf}
%Consider the Gelfand-Naimark semi norm on $A\otimes B$ defined as
%$$ \gamma(x) = \sup\{\|T(x)\|\},$$  i
%where the supremum runs over all the $*$-representations $T$ of   $A\otimes B$ on Hilbert spaces. Since the $\|\cdot\|_{\max}$-norm on  $A\otimes B$ is continuous with respect to $`\wedge'$-norm, it can be extended to $A\oop B$, and thus $A\otimes B$ is $\|\cdot\|_{\max}$-dense in $A\oop B$. Thus, by \cite[Proposition 10.5.20]{palmer}, $A\oop B$ is $*$-regular, whenever $A\otimes B$ is so. The particular case follows from \cite[Corollary 2.7]{hkv}.
%\end{pf}

\begin{rem}
 Let $A$ and $B$ be $C^*$-algebras. If $A$ is subhomogeneous, then by \cite[Proposition IV.1.4.6]{blackadar}, every subhomogeneous $C^*$-algebra is bidual type I, so that $A^{**}$ is a type $I$ von Neumann algebra, which is then nuclear by Corollary IV.2.2.10 and Theorem IV.3.1.5 of \cite{blackadar}. Consider the Gelfand-Naimark semi norm on $A\otimes B$ defined as
$$ \gamma(x) = \sup\{\|T(x)\|\},$$ 
where the supremum runs over all the $*$-representations $T$ of   $A\otimes B$ on Hilbert spaces. 
Since the $\|\cdot\|_{\gamma}$-norm on  $A\otimes B$ is continuous with respect to $`\wedge'$-norm, it can be extended to $A\oop B$, and thus $A\otimes B$ is $\|\cdot\|_{\gamma}$-dense in ($A\oop B, \|\cdot\|_{\gamma})$. By \cite[Proposition 10.5.20]{palmer}, $A\oop B$ is $*$-regular, whenever $A\otimes B$ is so. Since $A$ is nuclear, by \cite[Corollary 2.7]{hkv}, $A\oop B$ is $*$-regular.
%So, by Proposition \ref{spec11}, $A\oop B$ is $*$-regular.
It is also Hermitian (\cite[Theorem 4.6]{ranj3} and is $*$-semisimple (follows from \cite[Theorem 4.1]{ranj3}). Hence, the definition of spectral synthesis in \cite{somer2} is equivalent to our definition in this case. 
\end{rem}

In the case of commutative separable $C^*$-algebras the ideals which are not singly generated fail to be spectral. A similar result also holds true in the non-commutative situation. The following can be proved exactly on the same lines of \cite[Theorem 6.12]{ass}.

\begin{prop}
 Let $A$ and $B$ be separable $C^*$-algebras, and  $J$ be a non-zero closed ideal of $A\oop B$. Then $J$ is singly generated if it is spectral.
\end{prop}

\section{Reverse Involution}

Let $A$ be a $C^*$-algebra. On the Banach algebra $A\ot A$ (with usual multiplication), define the involution as $(a\ot b)^* = b^* \ot a^*$ for all $a,b\in A$. Then it extends to an isometric involution on $A\oop A$ and $A\oop A$ forms a Banach $*$-algebra with this involution, which we denote by $A\oop_r A$. Regarding the closed $*$-ideals of $A\oop_r A$, note that the closed ideals of $A\oop_r A$ coincide with the ones in $A\oop A$; however, the closed $*$-ideals differ. We do not know whether a closed ideal of $A\oop A$ is a $*$-ideal or not, but in $A\oop_r A$ a closed ideal need not be a $*$-ideal. For example, in the space $\B(H)\oop_r \B(H)$, the closed ideals $\K(H) \oop \B(H)$ and $\B(H) \oop \K(H)$ are not $*$-ideals. In fact, it has only two non-trivial closed $*$-ideals, namely  $\K(H) \oop_r \K(H)$ and $\B(H)\oop \K(H) + \K(H) \oop \B(H)$.

We know that with natural involution $A\oop A$ has a faithful $*$-representation and is always $*$-semisimple for any $C^*$-algebra $A$. However, we show that this is not the case with $A\oop_r A$. 

\begin{prop}
Let $A$ be a unital $C^*$-algebra. Then, $A\oop_r A$ has a faithful $*$-representation if and only if $A = \C I$, $I$ being the unity of $A$. 
\end{prop}

\begin{pf}
 Let $\pi$ be a faithful $*$-representation of $A\oop_r A$ on a Hilbert space $H$. Define $\pi_1(a):= \pi(1\ot a)$ and $\pi_2(a):=\pi(a\ot 1)$ for all $a\in A$. Then $\pi_1$ and $\pi_2$ are both bounded representations of $A$ on $\B(H)$, with $\pi(a\ot b) = \pi_1(a) \pi_2(b) = \pi_2(b) \pi_1(a)$ for all $a,b\in A$. Also 
\begin{equation}\label{reverse1}
 \pi_1(a^*) = \pi(1\ot a^*) = \pi ((a\ot 1)^*) = (\pi(a\ot 1))^* = \pi_2(a)^*
\end{equation}
for all $a\in A$. It is known that an element $h\in A$ is self adjoint if and only if $\|\exp ith\| =1$ for all $t \in \R$. For a self adjoint element $h\in A$, using the facts that $\pi $ is contractive
% being a $*$-representation of Banach $*$-algebra into a $C^*$-algebra, 
and that $\|\cdot\|_\wedge$-norm is a cross norm, we have
\begin{eqnarray*}
 \|\exp it \pi_1(h)\| &=& \|\pi (\exp it (h \ot 1))\| \\
               &\leq& \|\exp it (h\ot 1)\|_\wedge\\
               & = & \lim_m \Big\| \sum_{n=1}^m \frac{i^nt^n(h^n\ot 1)}{n!} \Big\|_\wedge\\
               & = & \lim_m \Big\| \Big(\sum_{n=1}^m \frac{i^nt^n h^n}{n!}\Big)\ot 1 \Big\|_\wedge\\
                & = & \lim_m \Big\| \sum_{n=1}^m \frac{i^nt^n h^n}{n!} \Big\|\\                
 & = & \|\exp ith\|=1, 
\end{eqnarray*}
and this is true for all $t\in \R$. Thus,  $\|\exp it \pi_1(h)\| =1$ for all $t\in \R$, which shows that $\pi_1(h)$ is a self adjoint element of $\B(H)$. This, combined with equation (\ref{reverse1}), gives $\pi_1(h) = \pi_2(h)$, that is $ \pi(1\ot h) = \pi(h\ot 1)$. Since $\pi$ is faithful, $1\ot h = h \ot 1$. So, for any $\phi \in A^*$,  $\phi(1) h = \phi(h) 1$, which further gives $h\in \C I$, and this is true for any self adjoint element $h$ of $A$. Since any $a\in A$ can be written as $a = h + ik$, $h$ and $k$ being self adjoint elements of $A$,  we obtain the required result.
\end{pf}

\begin{cor} 
\begin{enumerate}[(i)]
 \item  $A\oop_r A$ is $*$-semisimple if and only if $A = \C I$.
 \item $A \oop_r A$ is symmetric if and only if $A = \C I$ 
\end{enumerate}
\end{cor}

\begin{pf} (i) Follows easily from the fact that a semisimple Banach $*$-algebra possesses a faithful $*$-representation \cite[Corollary 4.7.16]{rick}. 

(ii) Let $A\oop_r A$ be symmetric. Using the same argument as in \cite[Proposition 5.16]{ass}, one can show that the radical of $A\oop_r A$ is $\{0\}$. By \cite[Theorem 4.7.15]{rick}, $*$-radical of $A\oop_r A$ coincides with its radical. Thus $A\oop_r A$ is $*$-semisimple, which using above part implies $A = \C I$.
\end{pf}

%====================================================================

\end{document}